\documentstyle[twoside,amstex]{article}
\input amssym.def
\input amssym.tex
\def\bc{\begin{center}}
\def\ec{\end{center}}
\def\no{\noindent}
\def\hang{\hangindent\parindent}
\def\textindent#1{\indent\llap{[#1]\enspace}\ignorespaces}
\def\re{\par\hang\textindent}
\begin{document}
\thispagestyle{empty} \vspace*{3 true cm} \pagestyle{myheadings}
\markboth {\hfill {\sl Huanyin Chen, H. Kose and Y.
Kurtulmaz}\hfill} {\hfill{\sl THE STRONG $P$-CLEANNESS OVER
RINGS}\hfill} \vspace*{-1.5 true cm} \bc{\large\bf THE STRONG
$P$-CLEANNESS OVER RINGS}\ec

\vskip6mm
\bc{{\bf Huanyin Chen}\\[1mm]
Department of Mathematics, Hangzhou Normal University\\
Hangzhou 310036, China, huanyinchen@@aliyun.com}\ec
\bc{{\bf H. Kose}\\[1mm]
Department of Mathematics, Ahi Evran University\\ Kirsehir,
Turkey, handankose@@gmail.com}\ec
\bc{{\bf Y. Kurtulmaz}\\[1mm]
Department of Mathematics, Bilkent University\\ Ankara, Turkey,
yosum@@fen.bilkent.edu.tr}\ec

\vskip10mm
\begin{minipage}{120mm}
\bc {\bf ABSTRACT}\ec

\vskip4mm An element of a ring $R$ is strongly $P$-clean provided
that it can be written as the sum of an idempotent and a strongly
nilpotent element that commute. A ring $R$ is strongly $P$-clean
in case each of its elements is strongly $P$-clean. We
investigate, in this article, the necessary and sufficient
conditions under which a ring $R$ is strongly $P$-clean. Many
characterizations of such rings are obtained. The criteria on
strong $P$-cleanness of $2\times 2$ matrices over commutative
local rings are also determined.\vskip3mm {\it Key Words:}\ \
Strongly $P$-clean ring; $n\times n$ Matrix; Local ring; Strongly
$\pi$-regularity.

\end{minipage}

\vskip20mm \bc{\bf 1. INTRODUCTION}\ec

\vskip4mm \no An element $a\in R$ is strongly clean provided that
there exist an idempotent $e\in R$ and an element $u\in U(R)$ such
that $a=e+u$ and $eu=ue$, where $U(R)$ is the set of all units in
$R$. A ring $R$ is strongly clean in case every element in $R$ is
strongly clean. Recently, strongly cleanness has been extensively
studied in the literature (cf. [1-5] and [10-13]). As is well
known, every $2\times 2$ matrix $A$ over a field satisfies the
conditions: $A=E+W, E$ is similar to a diagonal matrix, $W\in
M_2(R)$ is nilpotent and $E$ and $W$ commutate. Such a
decomposition over a field is called the Jordan-Chevalley
decomposition in Lie algebra theory. This motivates us to
investigate certain strong cleanness related to nilpotent
property. Following Diesl [8], a ring $R$ is strongly nil clean
provided that for any $a\in R$ there exists an idempotent $e\in R$
such that $a-e\in R$ is nilpotent and $ae=ea$. If such idempotent
is unique, we say $R$ is uniquely nil clean. In [6], the author
develop the theory for strongly nil clean matrices. The main
purpose of this article is to introduce a subclass of strongly nil
cleanness but behaving better than those ones.

An element $a$ of a ring $R$ is strongly nilpotent if every
sequence $a=a_0,a_1,a_2,\cdots $ such that $a_{i+1}\in a_iRa_i$ is
ultimately zero. Obviously, every strongly nilpotent element is
nilpotent. The prime radical $P(R)$ of a ring $R$, i.e. the
intersection of all prime ideals, consists of precisely the
strongly nilpotent elements. Replacing nilpotent elements by
strongly nilpotent elements, we shall investigate strong
$P$-cleanness over a ring $R$. An element of a ring $R$ is called
strongly $P$-clean provided that it can be written as the sum of
an idempotent and an element in $P(R)$ that commute. A ring $R$ is
strongly $P$-clean in case each of its elements is strongly
$P$-clean. In Section 2, we give several necessary and sufficient
conditions under which a ring $R$ is strongly $P$-clean. Many
characterizations of such rings are obtained. A ring $R$ is said
to be local if $R$ has only one maximal right ideal. In Section 3,
the strong $P$-cleanness of triangular matrix ring over a local
ring is determined. In Section 4, we characterize strongly
$P$-clean matrix over commutative local rings by means of the
solvability of quadratic equations. Finally, in Section 5, we
continuous to investigate such matrices via the characteristic
criteria.

Throughout, all rings are associative rings with identity. As
usual, $M_n(R)$ denotes the ring of all $n\times n$ matrices over
a ring $R$ and $GL_2(R)$ denotes the $2$-dimensional general
linear group of a ring $R$. An ideal $I$ of a ring $R$ is locally
nilpotent provided that for any $x\in I$, $RxR$ is nilpotent. Let
$a\in R$. Then $ann_{\ell}(a)=\{ r\in R~|~ra=0\}$ and
$ann_{r}(a)=\{ r\in R~|~ar=0\}$. $J(R)$ and $P(R)$ stand for the
Jacobson radical and prime radical of $R$, respectively.

\vskip10mm \bc{\bf 2. STRONGLY $P$-CLEAN RINGS}\ec

\vskip4mm \no Recall that a ring $R$ is Boolean provided that
every element in $R$ is an idempotent. Obviously, all Boolean
rings are commutative. Let $R$ be a ring. Then $P(R)=\{ x\in
R~|~RxR~\mbox{is nilpotent}\}.$ We begin with the connection
between strong $P$-cleanness and strong cleanness.

\vskip4mm \no{\bf Theorem 2.1.}\ \ {\it A ring $R$ is strongly
$P$-clean if and only if}\vspace{-.5mm}
\begin{enumerate}
\item [(1)]{\it $R$ is strongly clean.}
\vspace{-.5mm}
\item [(2)]{\it $R/J(R)$ is Boolean.}
\vspace{-.5mm}
\item [(3)]{\it $J(R)$ is locally nilpotent.}
\vspace{-.5mm}
\end{enumerate}\no{\it Proof.}\ \ Suppose that $R$ is strongly $P$-clean. Let $x\in R$. Then
there exists an idempotent $e\in R$ and a $w\in P(R)$ such that
$x=e+w$ and $ew=we$. Thus, $x=1+(w-1)$. Write $w^m=0$. Then
$(w-1)^{-1}=-1-w-\cdots -w^{m-1}$. Hence, $x\in R$ is strongly
clean. Thus,  $R$ is strongly clean. Clearly, $P(R)\subseteq
J(R)$. This implies that $R/J(R)$ is Boolean. Let $x\in J(R)$.
Then there exists an idempotent $e\in R$ and an element $w\in
P(R)$ such that $x=e+w$. Clearly, $w\in J(R)$, and so $e=x-w\in
J(R)$. This implies that $e=0$. Hence, $x=w\in P(R)$, i.e., $RxR$
is nilpotent. Therefore $J(R)$ is locally nilpotent.

Conversely, assume that conditions $(1), (2)$ and $(3)$ hold. Let
$x\in R$. Since $R$ is strongly clean, we can find an idempotent
$e\in R$ and an invertible $u\in R$ such that $x=e+u$ and $ex=xe$.
Thus, $x=(1-e)+(2e-1+u)$ and $(1-e)^2=1-e$. As $R/J(R)$ is
Boolean, we see that $\overline{u}^2=\overline{u}$, and so $u-1\in
J(R)$. As $\overline{2}^2=\overline{2}\in R/J(R)$, we deduce that
$2\in J(R)$; hence, $2e-1+u\in J(R)$. Since $J(R)$ is locally
nilpotent, $R(2e-1+u)R$ is nilpotent; hence, $2e-1+u\in P(R)$, as
required.\hfill$\Box$

\vskip4mm Recall that a ring $R$ is strongly $J$-clean provided
that for any $x\in R$, there exists an idempotent $e\in R$ such
that $x-e\in J(R)$ and $xe=ex$ (cf. [5]). One easily checks that a
ring $R$ is strongly $P$-clean if and only if $R$ is strongly
$J$-clean and $J(R)$ is locally nilpotent.

\vskip4mm \no{\bf Corollary 2.2.}\ \ {\it Let $R$ be a local ring.
Then the following are equivalent:}\vspace{-.5mm}
\begin{enumerate}
\item [(1)]{\it $R$ is strongly $P$-clean.}
\vspace{-.5mm}
\item [(2)]{\it $R/J(R)\cong {\Bbb Z}_2$ and $J(R)$ is locally nilpotent.}
\end{enumerate}\no{\it Proof.}\ \ It is immediate from Theorem 2.1.\hfill$\Box$

\vskip4mm The following example shows that strongly clean may be
not strongly $P$-clean.

\vskip4mm \no{\bf Example 2.3.}\ \ Let
$R=\prod\limits_{n=1}^{\infty}{\Bbb Z}_{2^n}$. For each $n$,
${\Bbb Z}_{2^n}$ is a local ring with the Jacobson radical $2{\Bbb
Z}_{2^n}$. One easily checks that ${\Bbb Z}_{2^n}$ is strongly
clean. Thus, $R$ is strongly clean. Choose $r=(0,2,2,2,\cdots )$.
It is easy to check that $r\in R$ is not strongly $P$-clean.
Therefore $R$ is not a strongly $P$-clean ring.

\vskip4mm Let $comm(x)=\{r\in R~|~xr=rx\}$ and $comm^2(x)=\{r\in
R~|~ry=yr~\mbox{for all}~y\in comm(r)\}$.

\vskip4mm \no{\bf Theorem 2.4.}\ \ {\it Let $R$ be a ring. Then
the following are equivalent}\vspace{-.5mm}
\begin{enumerate}
\item [(1)]{\it $R$ is strongly $P$-clean.}
\vspace{-.5mm}
\item [(2)]{\it $R/P(R)$ is Boolean.}
\vspace{-.5mm}
\item [(3)]{\it For any $x\in R$, there exists an idempotent $e\in R$ such
that $x-e\in P(R)$.} \vspace{-.5mm}
\item [(4)]{\it For any $x\in R$, there exists an idempotent $e\in comm^2(x)$ such
that $x-e\in P(R)$.} \vspace{-.5mm}
\item [(5)]{\it For any $x\in R$, there exists a unique idempotent $e\in R$ such
that $x-e\in P(R)$.} \vspace{-.5mm}
\end{enumerate}\no{\it Proof.}\ \ $(1)\Rightarrow (3)$ is trivial.

$(3)\Rightarrow (2)$ is clear.

$(2)\Rightarrow (4)$ By hypothesis, $R/P(R)$ is Boolean. For any
$x\in R$, then $\overline{x}\in R/P(R)$ is an idempotent. Hence,
$x-x^2\in P(R)$, i.e., $x(1-x)\in P(R)$. Write $x^n(1-x)^n=0$. Let
$f(t)=\sum\limits_{i=0}^n \left(
\begin{array}{c}
2n\\
i
\end{array}
\right)t^{2n-i}(1-t)^i\in {\Bbb Z}[t]$. Then $f(t)\equiv 0
\,(mod\,t^n)$. It follows from $$f(t)+ \sum\limits_{i=n+1}^{2n}
\left( \begin{array}{c} 2n\\
i \end{array} \right)x^{2n-i}(1-t)^i=\big(t+(1-t)\big)^n=1$$ that
$f(t)\equiv 1 \,\big(mod\,(1-t)^n\big)$. Thus,
$f(t)\big(1-f(t)\big)\equiv 0\,\big(mod\,t^n(1-t)^n\big)$. Let
$e=f(x)$. We see that $e(1-e)=0$; hence, $e\in R$ is an
idempotent. For any $y\in comm(x)$, we have $yx=xy$, and then
$ye=yf(x)=f(x)y=ey$. This implies that $y\in comm^2(x)$. Further,
$x-e\in P(R)$, and therefore $R$ is strongly $P$-clean.

$(4)\Rightarrow (1)$ As $x\in comm(x)$, this implication is obvious.

$(5)\Leftrightarrow (1)$ This is easily obtained.\hfill$\Box$

\vskip4mm Immediately, we see that every Boolean ring is strongly
$P$-clean. As every Boolean ring has stable range one, it follows
from Theorem 2.4 that every strongly $P$-clean ring has stable
range one.

\vskip4mm \no{\bf Corollary 2.5.}\ \ {\it A ring $R$ is strongly
$P$-clean if and only if}\vspace{-.5mm}
\begin{enumerate}
\item [(1)]{\it $R$ is periodic;}
\vspace{-.5mm}
\item [(2)]{\it Every element in $1+U(R)$ is strongly nilpotent.}
\end{enumerate}\no{\it Proof.}\ \ Suppose $R$ is strongly
$P$-clean. For any $x\in R$, it follows by Theorem 2.4 that
$x-x^2\in P(R)$. Thus, $(x-x^2)^n=0$ for some $n\in {\Bbb N}$.
This shows that $x^n=x^{n+1}f(x)$, where $f(t)\in {\Bbb Z}[t]$. By
using Herstein's Theorem, $R$ is periodic. Let $x\in 1+U(R)$.
Write $x=e+w$ with $e=e^2, w\in P(R)$ and $ae=ea$. Then
$1-x=(1-e)-w$, and so $1-e=(1-x)+w\in U(R)$. It follows that
$e=0$, and therefore $x=w\in P(R)$ is strongly nilpotent.

Conversely, assume that $(1)$ and $(2)$ hold. Since $R$ is
periodic, it is strongly $\pi$-regular. In view of [5, Proposition
13.1.8], there exist $e=e^2\in R, u\in U(R)$ and a nilpotent $w\in
R$ such that $x=eu+w$, where $e,u,w$ commutate. By hypothesis,
$1-u\in P(R)$, and then $u\in 1+P(R)$. Moreover, we see that
$w=1-(1-w)\in P(R)$. Accordingly, $x=e+\big(w-x(1-u)big)$ with
$w-x(1-u)\in P(R)$. Therefore $R$ is strongly $P$-clean.
\hfill$\Box$

\vskip4mm Let ${\Bbb Z}_{2^n}[i]=\{ a+bi~|~a,b\in {\Bbb
Z}_{2^n},i^2=-1\} (n\geq 2)$. Then we claim that ${\Bbb
Z}_{2^n}[i]$ is strongly $P$-clean. One easily checks that
$P({\Bbb Z}_{2^n}[i])=(1+i)$. Further, ${\Bbb
Z}_{2^n}[i]/P\big({\Bbb Z}_{2^n}[i]\big)\cong {\Bbb Z}_2$ is
Boolean, and we are through by Theorem 2.4.

Let $R= \left(
\begin{array}{cc}
{\Bbb Z}_2&{\Bbb
Z}_2\\
0&{\Bbb Z}_2 \end{array} \right)$. Then $P(R)=\left(
\begin{array}{cc}
0&{\Bbb Z}_2\\
0&0 \end{array} \right)$. Hence, $R/P(R)\cong {\Bbb Z}_2\oplus
{\Bbb Z}_2$, and so $R/P(R)$ is Boolean. Therefore $R$ is strongly
$P$-clean. Let $R=\{ \left(
\begin{array}{cc}
a&b\\
0&a \end{array} \right)~|~a\in {\Bbb Z}_2, b\in {\Bbb Z}\}$. Then
$R/P(R)\cong {\Bbb Z}_2$ where $P(R)=\left(
\begin{array}{cc}
0&{\Bbb Z}\\
0&0 \end{array} \right)$. This implies that $R$ is strongly
$P$-clean.

\vskip4mm \no{\bf Lemma 2.6.}\ \ {\it Every homomorphic image of
strongly $P$-clean rings is strongly $P$-clean.} \vskip2mm\no{\it
Proof.}\ \ Let $I$ be an ideal of a strongly $P$-clean ring $R$.
Let $M$ be a prime ideal of $R/I$. Then $M=P/I$, where $P$ is a
prime ideal of $R$. Let $\overline{x}\in R/I$. In light of Theorem
2.4, $x-x^2\in P$; hence, $\overline{x}-\overline{x}^2\in M$. This
shows that $\overline{x}-\overline{x}^2\in P\big(R/I\big)$. Thus
$R/I/P\big(R/I\big)$ is Boolean, and we therefore complete the
proof by Theorem 2.4.\hfill$\Box$

\vskip4mm \no{\bf Lemma 2.7.}\ \ {\it Let $I$ be a nilpotent ideal
of a ring $R$. Then $R$ is strongly $P$-clean if and only if $R/I$
is strongly $P$-clean.} \vskip2mm\no{\it Proof.}\ \ If $R$ is
strongly $P$-clean, then so is $R/I$ by Lemma 2.6. Write $I^n=0
(n\in {\Bbb N})$. Suppose $R/I$ is strongly $P$-clean. For any
$x\in R$, it suffice to show that $x-x^2\in P(R)$ by Theorem 2.4.
Given $x-x^2=a_0,a_1,\cdots ,a_n,\cdots $ with each $a_{i+1}\in
a_iRa_i$, we have
$\overline{x-x^2}=\overline{a_0},\overline{a_1},\cdots
,\overline{a_n},\cdots $ with each $\overline{a_{i+1}}\in
\overline{a_i}(R/I)\overline{a_i}$. As $R/I$ is strongly
$P$-clean, it follows by Theorem 2.4 that
$\overline{a_m}=\overline{0}$ for some $m\in {\Bbb N}$. Hence,
$a_m\in I$. This shows that $a_{n+m}\in
\underbrace{\big(a_mR\big)\big(a_mR\big)\cdots
\big(a_mR\big)}\limits_{n}\subseteq I^n=0$. Therefore $x-x^2\in
P(R)$, hence the result. \hfill$\Box$

\vskip4mm \no{\bf Theorem 2.8.}\ \ {\it Let $I$ be an ideal of a
ring $R$. Then the following are equivalent:}\vspace{-.5mm}
\begin{enumerate}
\item [(1)]{\it $R/I$ is strongly $P$-clean.}
\vspace{-.5mm}
\item [(2)]{\it $R/I^n$ is strongly $P$-clean for some $n\in {\Bbb N}$.}
\vspace{-.5mm}
\item [(3)]{\it $R/I^n$ is strongly $P$-clean for all $n\in {\Bbb N}$.}
\end{enumerate}\no{\it Proof.}\ \ $(1)\Rightarrow (3)$ It is easy
to verify that $$R/I\cong \big(R/I^n\big)/\big(I/I^n\big).$$ As
$\big(I/I^n\big)^n=0$, we see that $R/I$ is strongly $P$-clean, by
Lemma 2.7.

$(3)\Rightarrow (2)$ is trivial.

$(2)\Rightarrow (1)$ Clearly,
$$R/I\cong \big(R/I^n\big)/\big(I/I^n\big).$$ Therefore the proof is complete in terms of
Lemma 2.6.\hfill$\Box$

\vskip4mm \no{\bf Lemma 2.9.}\ \ {\it Every finite subdirect of
strongly $P$-clean rings is strongly $P$-clean.} \vskip2mm\no{\it
Proof.}\ \ Let $R$ be the subdirect product of $R_1,\cdots ,R_n$,
where each $R_i$ is strongly $P$-clean. Then
$\bigoplus\limits_{i=1}^{n}R_i$ is strongly $P$-clean.
Furthermore, $R$ is a subring of $\bigoplus\limits_{i=1}^{n}R_i$.
Let $x\in R$. Then $x-x^2\in
P\big(\bigoplus\limits_{i=1}^{n}R_i\big)$. Given $x-x^2=a_0,
a_1,\cdots ,a_m,\cdots $ in $R$ and each $a_{i+1}\in a_iRa_i$, we
see that $x-x^2=a_0, a_1,\cdots ,a_m,\cdots $ in
$\bigoplus\limits_{i=1}^{n}R_i$ and each $a_{i+1}\in
a_i\big(\bigoplus\limits_{i=1}^{n}R_i\big)a_i$. In view of Theorem
2.4, $x-x^2\in P\big(\bigoplus\limits_{i=1}^{n}R_i\big)$. Hence,
we can find some $s\in {\Bbb N}$ such that $a_s=0$. This implies
that $x-x^2\in P(R)$. That is, $R/P(R)$ is Boolean. In light of
Theorem 2.4, $R$ is strongly $P$-clean, as required.\hfill$\Box$

\vskip4mm \no{\bf Proposition 2.10.}\ \ {\it Let $I$ and $J$ be
ideals of a ring $R$. Then the following are
equivalent:}\vspace{-.5mm}
\begin{enumerate}
\item [(1)]{\it $R/I$ and $R/J$ are strongly $P$-clean.}
\vspace{-.5mm}
\item [(2)]{\it $R/\big(IJ\big)$ is strongly $P$-clean.}
\vspace{-.5mm}
\item [(3)]{\it $R/\big(I\bigcap J\big)$ is strongly $P$-clean.}
\end{enumerate}\no{\it Proof.}\ \ $(1)\Rightarrow (3)$ Construct
maps $f: R/\big(I\bigcap J\big)\to R/I, x+\big(I\bigcap
J\big)\mapsto x+I$ and $g: R/\big(I\bigcap J\big)\to R/I,
x+\big(I\bigcap J\big)\mapsto x+J$. Then $ker(f)\bigcap ker(g)=0$.
Therefore $R/\big(I\bigcap J\big)$ is the subdirect product of
$R/I$ and $R/J$. Thus, $R/\big(I\bigcap J\big)$ is strongly
$P$-clean, by Lemma 2.9.

$(3)\Rightarrow (2)$ Obviously, $R/\big(I\bigcap J\big)\cong
\big(R/IJ\big)/\big((I\bigcap J)/IJ\big)$, and $\big((I\bigcap
J)/IJ\big)^2=0$. In view of Lemma 2.7, $R/\big(IJ\big)$ is
strongly $P$-clean.

$(2)\Rightarrow (1)$ As $R/I\cong \big(R/IJ\big)/\big(I/IJ\big)$,
it follows from Lemma 2.6 that $R/I$ is strongly $P$-clean.
Likewise, $R/J$ is strongly $P$-clean.\hfill$\Box$

\vskip4mm We say that a ring $R$ is uniquely $P$-clean provided
that for any $x\in R$ there exists a unique idempotent $e\in R$
such that $x-e\in P(R)$.

\vskip4mm \no{\bf Theorem 2.10.}\ \ {\it Let $R$ be a ring. Then
$R$ is uniquely $P$-clean if and only if}\vspace{-.5mm}
\begin{enumerate}
\item [(1)]{\it $R$ is abelian;}
\vspace{-.5mm}
\item [(2)]{\it $R$ is strongly $P$-clean.}
\vspace{-.5mm}
\end{enumerate}\no{\it Proof.}\ \ Suppose $R$ is uniquely $P$-clean.
For all $x\in R$ there exists a unique idempotent $e\in R$ such
that $x-e\in P(R)$. Thus, $R/P(R)$ is Boolean. In view of Theorem
2.4, $R$ is strongly $P$-clean. Furthermore,
$\overline{ex-exe}^2=\overline{ex-exe}=0$. Hence, $ex-exe\in
P(R)$. Clearly, $e$ and $e+ex-exe\in R$ are idempotents, and that
$e-e,e-(e+ex-exe)\in P(R)$. By the uniqueness, we get $ex=exe$.
Likewise, $xe=exe$, and so $ex=xe$. That is, every idempotent in
$R$ is central. Therefore $R$ is abelian.

Conversely, assume that $(1)$ and $(2)$ hold. For any $x\in R$,
there exists an idempotent $e\in R$ such that $x-e\in P(R)$.
Suppose that $x-f\in P(R)$ where $f\in R$ is an idempotent. Then
$e-f=(x-f)-(x-e)\in P(R)$. Hence, we can find some $n\in {\Bbb N}$
such that $(e-f)^{2n+1}=e-f=0$. This implies that $e=f$, as
required.\hfill$\Box$

\vskip4mm In light of Theorem 2.10,one directly verifies that
${\Bbb Z}_4$ is uniquely $P$-clean. Recall that a ring $R$ is a
uniquely clean ring provided that each element in $R$ has a unique
representation as the sum of an idempotent and a unit (cf. [12]).
Let $R= \left(
\begin{array}{cc}
{\Bbb Z}_2&{\Bbb Z}_2\\
0&{\Bbb Z}_2
\end{array}
\right)$. By [12, Example 21], $R$ is not uniquely clean. But it
is strongly $P$-clean.

\vskip4mm \no{\bf Corollary 2.11.}\ \ {\it Every uniquely
$P$-clean is uniquely clean.}\vskip2mm\no{\it Proof.}\ \ In view
of Theorem 2.1, $R$ is strongly clean. Write $x=e+u$ where
$e=e^2\in R$ and $u\in U(R)$. Then $(1-e)-x=(1-2e)-u$. Clearly,
$(1-2e)^2=1$. As $R/P(R)$ is Boolean, we see that
$\overline{u}=\overline{1-2e}=\overline{1}$. Thus, $(1-2e)-u\in
P(R)$. This implies that $(1-e)-x\in P(R)$. Write $x=f+v$ where
$f=f^2\in R$ and $v\in U(R)$. Likewise, $(1-f)-x\in P(R)$. By the
uniqueness, we get $1-e=1-e$, and then $e=f$. Therefore $R$ is
uniquely clean.\hfill$\Box$

\vskip4mm \no {\bf Corollary 2.12.}\ \ {\it Let $R$ be uniquely
$P$-clean. Then $T=\{ (a_{ij})\in T_n(R)~|~a_{11}=\cdots
=a_{nn}\}$ is strongly $P$-clean.} \vskip2mm\no {\it Proof.}\ \
Let $S=\{ (a_{ij})\in T_n(R)~|~a_{11}=\cdots =a_{nn}=0\}$. Then
$S$ be a ring (not necessary unitary), and $S$ is a
$R$-$R$-bimodule in which $(s_1s_2)r=s_1(s_2r),
r(s_1s_2)=(rs_1)s_2$ and $(s_1r)s_2=s_1(rs_2)$ for all $s_1,s_2\in
S, r\in R$. Construct $I(R;S)=\{ (r,s)~|~r\in R,s\in S\}$. Define
$(r_1,s_1)+(r_2,s_2)=(r_1+r_2,s_1+s_2);
(r_1,s_1)(r_2,s_2)=(r_1r_2,s_1s_2+r_1s_2+s_1r_2).$ Then $I(R;S)$
is a ring with an identity $(1,0)$. Obviously, $T\cong I(R;S)$.
Let $(r,s)\in I(R;S)$. Since $R$ is strongly $P$-clean, write
$r=e+w, ew=we, e=e^2\in R, w\in P(R)$. Hence, $(r,s)=(e,0)+(w,s)$.
Clearly, $(e,0)^2=(e,0)$. In light of Theorem 2.10, every
idempotent in  $R$ is central, we see that $es=se$, and so
$(e,0)(w,s)=(w,s)(e,0)$. As $w\in P(R)$, we can find some $m\in
{\Bbb N}$ such that $(RwR)^m=0$. This implies that
$\big(I(R;S)(w,s)I(R;S)\big)^{m+n}=(0,0)$. Hence, $(w,s)\in
P\big(I(R;S)\big)$. Therefore $I(R;S)$ is strongly $P$-clean, as
required.\hfill$\Box$

\vskip4mm \no{\bf Theorem 2.13.}\ \ {\it Let $R$ be a ring. Then
$R$ is uniquely $P$-clean if and only if}\vspace{-.5mm}
\begin{enumerate}
\item [(1)]{\it $R$ is strongly $P$-clean.}
\vspace{-.5mm}
\item [(2)]{\it $R$ is uniquely nil clean.}
\vspace{-.5mm}
\end{enumerate}\no{\it Proof.}\ \ Suppose $R$ is uniquely
$P$-clean. It follows by Theorem 2.10 that $R$ is strongly
$P$-clean. Additionally, $R$ is abelian. Let $w\in R$ is
nilpotent. Then we have an idempotent $e\in R$ such that $w-e\in
P(R)$ and $we=ew$. This shows that $e=w-(w-e)\in R$ is nilpotent.
Hence, $e=0$, and so $w\in P(R)$. Therefore $R$ is uniquely nil
clean.

Conversely, assume that $(1)$ and $(2)$ hold. Then $R$ is abelian.
Therefore we complete the proof by Theorem 2.10. \hfill$\Box$

\vskip4mm We note that $\{$~$R$ is uniquely $P$-
clean~$\}\subsetneqq $ $\{$~strongly $P$-clean rings
$\}\subsetneqq \{$ strongly clean rings $\}$.

\vskip4mm \no{\bf Corollary 2.14.}\ \ {\it Let $R$ be a ring. Then
$R$ is Boolean if and only if}\vspace{-.5mm}
\begin{enumerate}
\item [(1)]{\it $R$ is uniquely $P$-clean.}
\vspace{-.5mm}
\item [(2)]{\it Every primary ideal of $R$ is prime.}
\vspace{-.5mm}
\end{enumerate}\no{\it Proof.}\ \ One direction is obvious.

Conversely, assume that $(1)$ and $(2)$ holds. In view of Theorem
2.13, $R$ is uniquely nil clean. As every primary ideal of $R$ is
prime, we get $P(R)=\bigcap \{ P~|~P~\mbox{is primary}~\}$.
Similarly to [7, Lemma 4.6], we see that $P(R)=0$. Therefore $R$
is Boolean in terms of Theorem 2.4. \hfill$\Box$

\vskip10mm \bc{\bf 3. TRIANGULAR MATRIX RINGS}\ec

\vskip4mm \no We use $T_n(R)$ to denote the ring of all upper
triangular $n\times n$ matrix over a ring $R$. The aim of this
section is to investigate the conditions under which $T_n(R)$ is
strongly $P$-clean for a local ring $R$.

\vskip4mm \no{\bf Lemma 3.1.}\ \ {\it Let $R$ be a ring, and let
$a=e+w$ be a strongly $P$-clean decomposition of $a$ in $R$. Then
$ann_{\ell}(a)\subseteq ann_{\ell}(e)$ and $ann_{r}(a)\subseteq
ann_{r}(e)$.}\vskip2mm\no {\it Proof.}\ \ Let $r\in
ann_{\ell}(a)$. Then $ra=0$. Write $a=e+w, e=e^2,w\in P(R)$ and
$ew=we$. Then $re=-rw$; hence, $re=-rwe=-rew$. It follows that
$re(1+w)=0$ as $1+w\in U(R)$, and so $re=0$. That is, $r\in
ann_{\ell}(e)$. Therefore $ann_{\ell}(a)\subseteq ann_{\ell}(e)$.
A similar argument shows that $ann_{r}(a)\subseteq
ann_{r}(e)$.\hfill$\Box$

\vskip4mm \no{\bf Theorem 3.2.}\ \ {\it Let $R$ be a ring, and let
$f\in R$ be an idempotent. Then $a\in fRf$ is strongly $P$-clean
in $R$ if and only if $a\in fRf$ is strongly $P$-clean in
$fRf$.}\vskip2mm\no {\it Proof.}\ \ Suppose that $a=e+w,e=e^2\in
fRf, w\in P(fRf)$ and $ew=we$. Then there exists some $n\in {\Bbb
N}$ such that $(fRfwfRf)^n=0$, and so $(RfwfR)^{n+4}=0$. That is,
$(RwR)^{n+4}=0$. This infers that $w\in P(R)$. Hence, $a\in fRf$
is strongly $P$-clean in $R$.

Conversely, suppose that $a=e+w,e=e^2\in R,w\in P(R)$ and $ew=we$.
As $a\in fRf$, it follows from Lemma 3.1 that $$\begin{array}{ccl}
1-f&\in
&ann_{\ell}(a)\bigcap ann_{r}(a)\\
&\subseteq &ann_{\ell}(e)\bigcap ann_{r}(e)\\
&=&R(1-e)\bigcap (1-e)R\\
&=&(1-e)R(1-e).
\end{array}$$ Hence, $ef=e=fe$. We observe that
$a=fef+fwf$, $(fef)^2=fef$. Furthermore, $fef\cdot
fwf=fewf=fwef=fwf\cdot fef$. As $w\in P(R)$, there exists some
$n\in {\Bbb N}$ such that $(RwR)^n=0$. Thus, $(fRfwfRf)^n\subseteq
(RwR)^n=0$, and so $fwf\in P(fRf)$. Therefore we complete the
proof.\hfill$\Box$

\vskip4mm As is well known, every corner of a strongly clean ring is
strongly clean. Analogously, we can derive
the following.

\vskip4mm \no{\bf Corollary 3.3.}\ \ {\it A ring $R$ is strongly
$P$-clean if and only if so is $eRe$ for all idempotents $e\in
R$.} \vskip2mm\no {\it Proof.}\ \ One direction is obvious. Let
$a\in eRe$. As $R$ is strongly $P$-clean, we see that $a\in eRe$
is strongly $P$-clean in $R$. According to Theorem 3.2, $a\in eRe$
is strongly $P$-clean in $eRe$.\hfill$\Box$

\vskip4mm Let $a\in R$. Then $l_a: R\to R$ and $r_a: R\to R$ denote,
respectively, the abelian group endomorphisms given by $l_a(r)=ar$
and $r_a(r)=ra$ for all $r\in R$. Thus, $l_a-r_b$ is an abelian
group endomorphism such that $(l_a-r_b)(r)=ar-rb$ for any $r\in R$.

\vskip4mm \no{\bf Lemma 3.4.}\ \ {\it Let $R$ be a local ring and
suppose that $A=(a_{ij})\in T_n(R)$. Then for any set $\{ e_{ii}\}$ of
idempotents in $R$ such that $e_{ii}=e_{jj}$ whenever
$l_{a_{ii}}-r_{a_{jj}}$ is not a surjective abelian group
endomorphism of $R$, there exists an idempotent $E\in T_n(R)$ such
that $AE=EA$ and $E_{ii}=e_{ii}$ for every $i\in \{ 1,\cdots
,n\}$.}\vskip2mm\no {\it Proof.}\ \ See [1, Lemma 7].\hfill$\Box$

\vskip4mm \no{\bf Theorem 3.5.}\ \ {\it Let $R$ be a local ring.
Then the following are equivalent:}\vspace{-.5mm}
\begin{enumerate}
\item [(1)]{\it $R$ is strongly $P$-clean.}
\vspace{-.5mm}
\item [(2)]{\it $R$ is uniquely $P$-clean.}
\vspace{-.5mm}
\item [(3)]{\it $R/J(R)\cong {\Bbb Z}_2$ and $J(R)$
is locally nilpotent.} \vspace{-.5mm}
\item [(4)]{\it $T_n(R)$ is strongly $P$-clean.}
\vspace{-.5mm}
\end{enumerate}\no{\it Proof.}\ \ $(1)\Rightarrow (2)$ is obvious
from Theorem 2.10.

$(2)\Rightarrow (3)$ In view of Theorem 2.1, $R/J(R)$ is Boolean,
and $J(R)$ is local nilpotent. As $R$ is local, we get
$R/J(R)\cong {\Bbb Z}_2$.

$(3)\Rightarrow (4)$ Let $A=(a_{ij})\in T_n(R)$. We need to
construct an idempotent $E\in T_n(R)$ such that $EA=AE$ and such
that $A-E\in P\big(T_n(R)\big)$. By hypothesis, $R/J(R)\cong {\Bbb
Z}_2$ and $J(R)$ is locally nilpotent. Thus, $R=J(R)\bigcup
\big(1+J(R)\big)$. Begin by constructing the main diagonal of $E$.
Set $e_{ii}=0$ if $a_{ii}\in J(R)$, and set $e_{ii}=1$ otherwise.
Thus, $a_{ii}-e_{ii}\in J(R)$ for every $i$. If $e_{ii}\neq
e_{jj}$, then it must be the case (without loss of generality)
that $a_{ii}\in U(R)$ and $a_{jj}\in J(R)$. Thus, $a_{jj}\in P(R)$
is nilpotent. Write $a_{jj}^m=0$. Construct a map
$\varphi=l_{a_{ii}^{-1}}+l_{a_{ii}^{-2}}r_{a_{jj}}+\cdots
+l_{a_{ii}^{-m}}r_{a_{jj}^{m-1}}: R\to R$. For any $r\in R$, it is
easy to verify that $\big(l_{a_{ii}}-r_{a_{jj}}\big)\big(\varphi
(r)\big)=r$. Thus, $l_{a_{ii}}-r_{a_{jj}}: R\to R$ is surjective.
According to Lemma 3.4, there exists an idempotent $E\in T_n(R)$
such that $AE=EA$ and $E_{ii}=e_{ii}$ for every $i\in \{ 1,\cdots
,n\}$. Further, $a_{ii}-e_{ii}\in P(R)$. Write
$\big(R(a_{ii}-e_{ii})R\big)^{m_i}=0$. Then one easily checks that
$$\big(T_n(R)(A-E)T_n(R)\big)^{\sum\limits_{i=1}^{n}m_i+n+1}=0.$$
This implies that $A-E\in P\big(T_n(R)\big)$. Therefore $T_n(R)$
is strongly $P$-clean.

$(4)\Rightarrow (1)$ is clear by Corollary 3.3.\hfill$\Box$

\vskip4mm \no{\bf Corollary 3.6.}\ \ {\it Let $R$ be a local ring.
Then the following are equivalent:}\vspace{-.5mm}
\begin{enumerate}
\item [(1)]{\it $R$ is strongly $P$-clean.}
\vspace{-.5mm}
\item [(2)]{\it For any $A\in
T_2(R)$, $A\in P\big(T_2(R)\big)$ or $I_2-A\in P\big(T_2(R)\big)$
or there exists $P\in U\big(T_2(R)\big)$ such that $P^{-1}AP=
\left(\begin{array}{cc}
\lambda&0\\
0&\mu
\end{array}
\right)$, where $\lambda\in 1+P(R),\mu\in P(R)$.} \vspace{-.5mm}
\end{enumerate}\no{\it Proof.}\ \ $(1)\Rightarrow (2)$ In view of Theorem
2.1, $J(R)$ is locally nilpotent; hence, $J(R)=P(R)$. Further,
$R/P(R)\cong {\Bbb Z}_2$. In view of Theorem 3.5, $T_2(R)$ is
uniquely $P$-clean. Let $A\in T_2(R)$. If $A, I_2-A\not\in
P\big(T_2(R)\big)$, without loss of generality, we may assume that
$A=\left(\begin{array}{cc}
a&v\\
0&b
\end{array}
\right)$, where $a\in 1+P(R),b\in P(R)$ and $v\in R$. As $R$ is
local, we can find an idempotent $E=\left(\begin{array}{cc}
1&w\\
0&0
\end{array}
\right)\in T_2(R)$ such that $A-E\in P\big(T_2(R)\big)$ and
$AE=EA$. Let $P=\left(\begin{array}{cc}
1&-w\\
0&1
\end{array}
\right)$. Then $P^{-1}EP=\left(\begin{array}{cc}
1&0\\
0&0
\end{array}
\right)$. Hence, $P^{-1}AP=\left(\begin{array}{cc}
1&0\\
0&0
\end{array}
\right)+W$, where $W=(w_{ij})\in P\big(T_2(R)\big)$. Obviously,
$\left(\begin{array}{cc}
1&0\\
0&0
\end{array}
\right)W=W\left(\begin{array}{cc}
1&0\\
0&0
\end{array}
\right)$, and so $w_{12}=0$. Thus, $P^{-1}AP=\left(
\begin{array}{cc}
1+w_{11}&0\\
0&w_{22}
\end{array}
\right)$ where $w_{11},w_{22}\in P(R)$, as desired. Conversely,
let $A\in T_2(R)$. It is easy to verify that $A$ is strongly
$P$-clean, and therefore $T_2(R)$ is strongly $P$-clean.
Accordingly, $R$ is strongly $P$-clean, by Theorem
3.5.\hfill$\Box$

\vskip4mm We close this section by considering a single $2\times
2$ strongly $P$-clean triangular matrix over a local ring.

\vskip4mm \no{\bf Proposition 3.7.}\ \ {\it Let $R$ be a local ring, let $A=\left(
\begin{array}{ll}
a&v\\
0&b
\end{array}
\right)\in T_2(R)$. Then $A$ is strongly $P$-clean if and only if
$a$ and $b$ are in $P(R)$ or $1+P(R)$.} \vskip2mm \no{\it Proof.}\
\ Suppose that $A$ is strongly $P$-clean and $A, I_2-A\not\in
P\big(T_2(R)\big)$. Then there exists some $E=\left(
\begin{array}{cc}
e&w\\
0&f
\end{array}
\right)\in R$ such that $$\left(
\begin{array}{cc}
a&v\\
0&b
\end{array}
\right)-E\in P\big(T_2(R)\big)~\mbox{and}~\left(
\begin{array}{cc}
a&v\\
0&b
\end{array}
\right)E=E\left(
\begin{array}{cc}
a&v\\
0&b
\end{array}
\right).$$ Since $A$ and $B$ are local rings, we see that $e=0,1$
and $f=0,1$. Thus, $E=\left(
\begin{array}{cc}
1&x\\
0&0
\end{array}
\right)$ or $E=\left(
\begin{array}{cc}
0&x\\
0&1
\end{array}
\right)$ where $x\in R$. This implies that $a\in P(R), b\in
1+P(R)$ or $a\in 1+P(R), b\in P(R)$, as desired.

Suppose that $a, b\in P(R)$ or $a, b\in 1_A+P(R)$, then $A\in
M_2(R)$ is strongly $P$-clean. Assume that $a\in 1+P(R), b\in
P(R)$. As $P(R)$ is locally nilpotent, we may write $b^m=0$.
Construct a map $\varphi =l_{a^{-1}}+l_{a^{-2}}r_{b}+\cdots
+l_{a^{-m}}r_{b^{m-1}}: R\to R$. Choose $x=\varphi (v)$. Then one
easily checks that $\big(l_{a}-r_{b}\big)\big(\varphi (v)\big)=v$.
Hence, $ax-xb=v$. Choose $E=\left(
\begin{array}{cc}
1&x\\
0&0
\end{array}
\right)$. Then $E=E^2, A-E\in P\big(T_2(R)\big)$ and $$AE=\left(
\begin{array}{cc}
a&ax\\
0&0
\end{array}
\right)=\left(
\begin{array}{cc}
a&v+xb\\
0&0
\end{array}
\right)=EA.$$ Assume that $a\in P(R), b\in 1+P(R)$. Analogously,
we can find an idempotent $E\in T_2(R)$ such that $AE=EA$ and
$A-E\in P\big(T_2(R)\big)$. Therefore $A\in T_2(R)$ is strongly
$P$-clean.\hfill$\Box$

\vskip4mm Let ${\Bbb Z}_{3^n}[\alpha]=\{ a+b\alpha~|~a,b\in {\Bbb
Z}_{3^n}, \alpha^2+\alpha+1=0\} (n\geq 1)$. Then $P\big({\Bbb
Z}_{3^n}[\alpha]\big)=\big(1-\alpha\big)$, i.e., the principal
generated by $1-\alpha\in {\Bbb Z}_{3^n}[\alpha]$. Therefore
${\Bbb Z}_{2^n}[\alpha]$ is local. Additionally, $T_2\big({\Bbb
Z}_{3^n}[\alpha]\big)$ is not strongly $P$-clean, by Theorem 3.5.
But, we see from Proposition 3.7 that $\left(
\begin{array}{ll}
x&z\\
0&y
\end{array}
\right)\in T_2\big({\Bbb Z}_{3^n}[\alpha]\big)$ is strongly
$P$-clean if and only if $x,y\in \big(1-\alpha\big)$ or
$1+\big(1-\alpha\big)$.

\vskip10mm \bc{\bf 4. STRONGLY $P$-CLEAN MATRICES}\ec

\vskip4mm \no The main purpose of this section is to investigate
the strong $P$-cleanness of a single matrix over commutative local
rings.

\vskip4mm \no{\bf Lemma 4.1.}\ \ {\it Let $R$ be a ring. Then
$P\big(M_n(R)\big)=M_n\big(P(R)\big)$.} \vskip2mm \no{\it Proof.}\
\ Let $A=(a_{ij})\in P\big(M_n(R)\big)$.  We can find some $m\in
{\Bbb N}$ such that $\big(M_n(R)AM_n(R)\big)^m=0$. For any
$s_i,r_i\in R (1\leq i\leq n)$, we have $$\begin{array}{c}
\sum\limits_{i=1}^{n}\left(
\begin{array}{cccc}
s_i&0&\cdots &0\\
0&0&\cdots &0\\
\vdots&\vdots&\ddots &\vdots\\
0&0&\cdots &0
\end{array}
\right)A\left(
\begin{array}{cccc}
r_i&0&\cdots &0\\
0&0&\cdots &0\\
\vdots&\vdots&\ddots &\vdots\\
0&0&\cdots &0
\end{array}
\right)=\left(
\begin{array}{cccc}
\sum\limits_{i=1}^{n}s_ia_{11}r_i&0&\cdots &0\\
0&0&\cdots &0\\
\vdots&\vdots&\ddots &\vdots\\
0&0&\cdots &0
\end{array}
\right). \end{array}$$ Hence,
$\big(\sum\limits_{i=1}^{n}s_ia_{11}r_i\big)^m=0$. This implies
that $\big(Ra_{11}R\big)^m=0$, and so $a_{11}\in P(R)$. Likewise,
every $a_{ij}\in P(R)$. Thus, $A\in M_n\big(P(R)\big)$; hence,
$P\big(M_n(R)\big)\subseteq M_n\big(P(R)\big)$.

Given any $A=(a_{ij})\in M_n\big(P(R)\big)$, then each $a_{ij}\in
P(R)$. Since $a_{11}\in P(R)$, we can find some $m\in {\Bbb N}$
such that $\big(Ra_{11}R\big)^m=0$. For any $(s^{t}_{ij}),
(r^{t}_{ij})\in M_n(R) (1\leq t\leq p)$, we have
$$\sum\limits_{t=1}^{p}(s^{t}_{ij})\left(
\begin{array}{cccc}
a_{11}&0&\cdots &0\\
0&0&\cdots &0\\
\vdots&\vdots&\ddots &\vdots\\
0&0&\cdots &0
\end{array}
\right)(r^{t}_{ij})\in M_n\big(Ra_{11}R\big),$$  and so
$$\big(\sum\limits_{t=1}^{p}(s^{t}_{ij})\left(
\begin{array}{cccc}
a_{11}&0&\cdots &0\\
0&0&\cdots &0\\
\vdots&\vdots&\ddots &\vdots\\
0&0&\cdots &0
\end{array}
\right)(r^{t}_{ij})\big)^m=0.$$ This implies that $$\big(M_n(R)\left(
\begin{array}{cccc}
a_{11}&0&\cdots &0\\
0&0&\cdots &0\\
\vdots&\vdots&\ddots &\vdots\\
0&0&\cdots &0
\end{array}
\right)M_n(R)\big)^m=0.$$ Hence, $$\left(
\begin{array}{cccc}
a_{11}&0&\cdots &0\\
0&0&\cdots &0\\
\vdots&\vdots&\ddots &\vdots\\
0&0&\cdots &0
\end{array}
\right)\in P\big(M_n(R)\big).$$ Likewise,
$$\left(\begin{array}{cccc}
0&a_{12}&\cdots &0\\
0&0&\cdots &0\\
\vdots&\vdots&\ddots &\vdots\\
0&0&\cdots &0
\end{array}
\right), \cdots ,\left(\begin{array}{cccc}
0&0&\cdots &0\\
0&0&\cdots &0\\
\vdots&\vdots&\ddots &\vdots\\
0&0&\cdots &a_{nn}
\end{array}
\right)\in P\big(M_n(R)\big).$$ Therefore we have
$$A=\left(\begin{array}{cccc}
a_{11}&0&\cdots &0\\
0&0&\cdots &0\\
\vdots&\vdots&\ddots &\vdots\\
0&0&\cdots &0
\end{array}
\right)+\cdots +\left(\begin{array}{cccc}
0&0&\cdots &0\\
0&0&\cdots &0\\
\vdots&\vdots&\ddots &\vdots\\
0&0&\cdots &a_{nn}
\end{array}
\right)\in P\big(M_n(R)\big).$$ Consequently,
$M_n\big(P(R)\big)\subseteq P\big(M_n(R)\big)$, as
desired.\hfill$\Box$

\vskip4mm \no{\bf Theorem 4.2.}\ \ {\it Let $R$ be a local ring.
Then $A\in M_2(R)$ is strongly $P$-clean if and only if $A\in
M_2\big(P(R)\big)$ or $I_2-A\in M_2\big(P(R)\big)$ or $A$ is
similar to a matrix $\left(
\begin{array}{cc}
\lambda&0\\
0&\mu
\end{array} \right)$, where $\lambda\in P(R), \mu\in 1+P(R)$. } \vskip2mm\no
{\it Proof.}\ \ If $A\in M_2\big(P(R)\big)$ or $I_2-A\in
M_2\big(P(R)\big)$, it follows by Lemma 4.1 that either $A$ or
$I_2-A$ is in $P\big(M_2(R)\big)$, and so $A$ is strongly
$P$-clean. For any $w_1,w_2\in P(R)$, we see that $\left(
\begin{array}{cc}
1+w_1&0\\
0&w_2\end{array} \right)=\left(
\begin{array}{cc}
1&0\\
0&0\end{array} \right)+ \left(
\begin{array}{cc}
w_1&0\\
0&w_2\end{array} \right)$. In light of Lemma 4.1, $\left(
\begin{array}{cc}
w_1&0\\
0&w_2\end{array} \right)\in M_2\big(P(R)\big)$. Thus, one
direction is clear.

Conversely, assume that $A\in M_2(R)$ is strongly $P$-clean, and
that $A, I_2-A\not\in M_2\big(P(R)\big)$. Then there exists an
idempotent $E\in M_2(R)$ and a $W\in P\big(M_2(R)\big)$ such that
$A=E+W$ with $EW=WE$. This implies that the idempotent $E\neq 0,
I_2$. In view of [5, Lemma 16.4.11], $E$ is similar to $\left(
\begin{array}{cc}
0&w_1\\
1&1+w_2\end{array} \right)$. As $E=E^2$, we deduce that $w_1=w_2=0$; hence, $E$ is similar to $\left(
\begin{array}{cc}
0&0\\
1&1\end{array} \right)$.
Obviously,
$\left(
\begin{array}{cc}
1&0\\
1&1\end{array} \right)\left(
\begin{array}{cc}
0&0\\
1&1\end{array} \right)\left(
\begin{array}{cc}
1&0\\
-1&1\end{array} \right)=\left(
\begin{array}{cc}
0&0\\
0&1\end{array} \right).$ Thus, we have an
$H\in GL_2(R)$ such that $HEH^{-1}= \left(
\begin{array}{cc}
1&0\\
0&0\end{array} \right)$. Thus, $HAH^{-1}=\left(
\begin{array}{cc}
1&0\\
0&0\end{array} \right)+HWH^{-1}$. Set $V=(v_{ij}):=HWH^{-1}$. It
follows from $EW=WE$ that $\left(
\begin{array}{cc}
1&0\\
0&0\end{array} \right)V=V\left(
\begin{array}{cc}
1&0\\
0&0\end{array} \right)$; hence, $v_{12}=v_{21}=0$ and
$v_{11},v_{22}\in P(R)$. Therefore $A$ is similar to $\left(
\begin{array}{cc}
1+v_{11}&0\\
0&v_{22}\end{array} \right)$, as desired.\hfill$\Box$

\vskip4mm \no{\bf Lemma 4.3.}\ \ {\it Let $R$ be a local ring, and
let $A\in M_2(R)$ be strongly $P$-clean. Then $A\in
M_2\big(P(R)\big)$ or $I_2-A\in M_2\big(P(R)\big)$ or $A$ is
similar to a matrix $\left(
\begin{array}{cc}
0&\lambda\\
1&\mu\end{array} \right)$, where $\lambda\in P(R),\mu\in
1+P(R)$.}\vskip2mm\no{\it Proof.}\ \ If $A, I_2-A\not\in
M_2\big(P(R)\big)$, it follows from Theorem 4.2 that there exists
a $P\in GL_2(R)$ such that $P^{-1}AP=\left(
\begin{array}{cc}
\alpha&0\\
0&\beta\end{array} \right)$, where $\alpha\in 1+P(R),\beta\in
P(R)$. One computes that
$$\begin{array}{lll}
&&[\alpha-\beta,1]B_{12}\big(-\alpha(\alpha-\beta)^{-1}\big)B_{21}(1)P^{-1}APB_{21}(-1)B_{12}\big(
\alpha(\alpha-\beta)^{-1}\big)[(\alpha-\beta)^{-1},1]\\
&=&\left(
\begin{array}{cc}
0&-(\alpha-\beta)\alpha(\alpha-\beta)^{-1}\beta\\
1&(\alpha-\beta)\alpha(\alpha-\beta)^{-1}+\beta
\end{array}
\right).\end{array}$$ Here, $[\xi, \eta] = diag(\xi, \eta)$ and $B_{ij}(\xi )= I_2 +\xi E_{ij}$ where
$E_{ij}$ is the matrix with $1$ on the place $(i, j)$ and $0$ on
other places. Let
$\lambda=-(\alpha-\beta)\alpha(\alpha-\beta)^{-1}\beta$ and
$\mu=(\alpha-\beta)\alpha(\alpha-\beta)^{-1}+\beta$. Therefore $A$
is similar to $\left(
\begin{array}{cc}
0&\lambda\\
1&\mu\end{array} \right)$, where $\lambda\in P(R), \mu\in 1+P(R)$.
\hfill$\Box$

\vskip4mm \no{\bf Theorem 4.4.}\ \ {\it Let $R$ be a commutative
local ring. Then the following are equivalent:} \vspace{-.5mm}
\begin{enumerate}
\item [(1)]{\it $A\in M_2(R)$ is strongly $P$-clean.}
\vspace{-.5mm}
\item [(2)]{\it $A-A^2\in M_2\big(P(R)\big)$.} \vspace{-.5mm}
\item [(3)]{\it $A\in M_2\big(P(R)\big)$ or $I_2-A\in M_2\big(P(R)\big)$ or the equation
$x^2-trA\cdot x+detA=0$ has a root in $P(R)$ and a root in
$1+P(R)$.} \vspace{-.5mm}
\end{enumerate}\no{\it Proof.}\ \ $(1)\Rightarrow (2)$ Write $A=E+W$ with $EW=WE, W\in
N\big(M_2(R)\big)$. Then $A-A^2=W-EW-WE-W^2\in N\big(M_2(R)\big)$.
Therefore, $A-A^2\in M_2\big(P(R)\big))$, by Lemma 4.1.

$(2)\Rightarrow (1)$ Since $A-A^2\in M_2\big(P(R)\big)$, we get
$A-A^2\in P\big(M_2(R)\big)$ by Lemma 4.1. As $P\big(M_2(R)\big)$
is locally nilpotent, we can find an idempotent $E\in M_2(R)$ such
that $A-E\in P\big(M_2(R)\big)$. Explicitly, $AE=EA$, as required.

$(1)\Rightarrow (3)$ Let $A\in M_2(R)$ be strongly $P$-clean and
$A, I_2-A\not\in M_2\big(P(R)\big)$. By virtue of Lemma 4.2, $A$
is similar to the matrix $\left(
\begin{array}{cc}
\lambda&0\\
0&\mu\end{array} \right)\in M_2(R)$, where $\lambda\in P(R),\mu\in
1+P(R)$. Thus, $x^2-trA\cdot
x+detA=det(xI_2-A)=(x-\lambda)(x-\mu)$, which has a root
$\lambda\in P(R)$ and a root $\mu\in 1+P(R)$.

$(3)\Rightarrow (1)$ Let $A\in M_2(R)$. If $A\in
M_2\big(P(R)\big)$ or $I_2-A\in M_2\big(P(R)\big)$, it follows
from Lemma 4.1 that $A\in M_2(R)$ is strongly $P$-clean.
Otherwise, it follows by the hypothesis that the equation
$x^2-trA\cdot x+detA=0$ has a root $x_1\in P(R)$ and a root
$x_2\in 1+P(R)$. Clearly, $x_1-x_2\in -1+P(R)\subseteq U(R)$. In
addition, $trA=x_1+x_2\in 1+P(R)$ and $detA=x_1x_2\in P(R)$. As
$detA\in P(R)$, $A\not\in U(R)$. It follows from
$det(I_2-A)=1-trA+detA\in P(R)$ that $I_2-A\not\in GL_2(R)$. In
light of [11, Lemma 4], there are some $\lambda\in J(R), \mu\in
1+J(R)$ such that $A$ is similar to $B= \left(
\begin{array}{cc}
0&\lambda\\
1&\mu
\end{array}
\right)$. Further, $x^2-trB\cdot
x+detB=det(xI_2-B)=det(xI_2-A)=x^2-trA\cdot x+detA$; and so
$x^2-trB\cdot x+detB=0$ has a root in $1+P(R)$ and a root in
$P(R)$. In view of Lemma 4.3, there exists a $P\in GL_2(R)$ such
that $P^{-1}BP= \left(
\begin{array}{cc}
\alpha_1&0\\
0&\alpha_2\end{array} \right)$ for some $\alpha_1\in 1+P(R),
\alpha_2\in P(R)$. By virtue of Lemma 4.1, $P^{-1}BP=\left(
\begin{array}{cc}
1&0\\
0&0\end{array} \right)+\left(
\begin{array}{cc}
\alpha_1-1&0\\
0&\alpha_2\end{array} \right)$ is a strongly $P$-clean expression.
Consequently, $A\in M_2(R)$ is strongly $P$-clean.\hfill$\Box$

\vskip4mm Let $R=\{ \frac{m}{n}\in {\Bbb Q}~|~2~\nmid~n\}$. Then $R$
is a commutative local ring. Choose $A= \left( \begin{array}{cc}
1&2\\
3&2 \end{array} \right)\in M_2(R)$. Clearly, $A, I_2-A\not\in
M_2\big(P(R)\big)$. Further, the equation $x^2-trA\cdot x+detA=0$
has a root $4$ and a root $-1$. But $4, -1\not\in P(R)$. Thus,
$A\in M_2(R)$ is not strongly $P$-clean from Theorem 4.5. But
$A\in M_2(R)$ is strongly clean by [4, Corollary 2.2]. It is worth
noting that every strongly $P$-clean $2\times 2$ matrix over
integral domains must be an idempotent by Theorem 4.4.

\vskip4mm \hspace{-1.8em} {\bf Corollary 4.5.}\ \ {\it Let $R$ be
a commutative local ring, and let $A\in M_2(R)$. Then the
following are equivalent:} \vspace{-.5mm}
\begin{enumerate}
\item [(1)]{\it $A\in M_2(R)$ is
strongly $P$-clean.} \vspace{-.5mm}
\item [(2)]{\it $A\in M_2\big(P(R)\big)$ or $I_2-A\in M_2\big(P(R)\big)$, or $trA\in 1+P(R)$ and the equation $x^2-x=-\frac{detA}{tr^2A}$ has a root in $P(R)$.}
\end{enumerate}\vspace{-.5mm} {\it Proof.}\ \ $(1)\Rightarrow (2)$ Assume that $A, I_2-A\not\in M_2\big(P(R)\big)$. By virtue of Theorem 4.2, $A$ is
similar to $\left(
\begin{array}{cc}
\lambda&0\\
0&\mu
\end{array} \right)$, where $\lambda\in P(R), \mu\in 1+P(R)$. Thus,
$trA=\lambda+\mu\in 1+P(R), detA=\lambda\mu$. Clearly,
$y^2-(\lambda+\mu)y+\lambda\mu=0$ has a root in $P(R)$. Thus, so
does the equation
$$(\lambda+\mu)^{-1}y^2-y=-(\lambda+\mu)^{-1}\lambda\mu.$$ Set
$x=(\lambda+\mu)^{-1}y$. Then
$$(\lambda+\mu)x^2-(\lambda+\mu)x=-(\lambda+\mu)^{-1}\lambda\mu.$$ We
infer that $x^2-x=-(\lambda+\mu)^{-2}\lambda\mu.$ Therefore
$x^2-x=-\frac{detA}{tr^2A}$ has a root in $P(R)$.

$(2)\Rightarrow (1)$ Suppose that $trA\in 1+P(R)$ and the equation
$x^2-x=-\frac{detA}{tr^2A}$ has a root $a\in P(R)$. Then
$\big((trA)a\big)^2-trA\cdot \big((trA)a\big)+detA=0$. Thus, the
equation $x^2-tr(A)x+det(A)=0$ has a root $tr(A)a\in P(R)$. One
easily checks that $tr(A)(1-a)\in 1+P(R)$ is an other root of the
preceding equation. Therefore, $A|in M_2(R)$ is strongly $P$-clean
in terms of Theorem 4.4.\hfill$\Box$

\vskip4mm \no{\bf Example 4.6.}\ \ {\it Let ${\Bbb Z}_4=\{
\overline{0}, \overline{1}, \overline{2}, \overline{3}\}$, and let
$A=\left( \begin{array}{cc}
\overline{1}&\overline{2}\\
\overline{2}&\overline{2} \end{array} \right)\in M_2({\Bbb Z}_4)$.
Then $A-A^2=\left( \begin{array}{cc}
\overline{0}&\overline{0}\\
\overline{0}&\overline{2} \end{array} \right)\in M_2\big(P({\Bbb
Z}_4)\big)$. Thus, $A\in M_2({\Bbb Z}_4)$ is strongly $P$-clean.
In fact, we have the strongly $P$-clean decomposition: $A=\left(
\begin{array}{cc}
\overline{1}&\overline{2}\\
\overline{2}&\overline{0} \end{array} \right)+\left( \begin{array}{cc}
\overline{0}&\overline{0}\\
\overline{0}&\overline{2} \end{array} \right)$. In this case, $A,
I_2-A\not\in P\big(M_2({\Bbb Z}_4)\big)$.}\hfill$\Box$

\vskip10mm \bc{\bf 5. CHARACTERISTIC CRITERIA}\ec

\vskip4mm \no For several kinds of $2\times 2$ matrices over commutative local rings, we can
derive accurate characterizations.

\vskip4mm \no {\bf Theorem 5.1.}\ \ {\it Let $R$ be a commutative
local ring, and let $A\in M_2(R)$. If $A$ is strongly $P$-clean,
then either $A\in M_2\big(P(R)\big)$, or $I_2-A\in
M_2\big(P(R)\big)$, or $trA\in 1+P(R)$ and $tr^2A-4detA=u^2$ for a
$u\in 1+P(R)$.} \vskip2mm\hspace{-1.8em} {\it Proof.}\ \ According
to Theorem 4.5, $A\in M_2\big(P(R)\big)$ or $I_2-A\in
M_2\big(P(R)\big)$, or $trA\in 1+P(R)$ and the equation
$x^2-x=\frac{detA}{-tr^2A}$ has a root $a\in P(R)$. Then $detA\in
P(R)$ and $2a-1\in -1+P(R)$. Further,
$(2a-1)^2=4(a^2-a)+1=\frac{4detA}{-tr^2A}+1=\frac{tr^2A-4detA}{tr^2A}$,
and therefore $tr^2A-4detA=\big(trA\cdot (2a-1)\big)^2$. Set
$u=trA\cdot (2a-1)$. Then $u\in 1+P(R)$, as required. \hfill$\Box$

\vskip4mm \no {\bf Corollary 5.2.}\ \ {\it Let $R$ be a commutative
local ring. If $\frac{1}{2}\in R$, then the following are equivalent:} \vspace{-.5mm}
\begin{enumerate}
\item [(1)]{\it $A\in M_2(R)$ is strongly $P$-clean.} \vspace{-.5mm}
\item [(2)]{\it $A\in M_2\big(P(R)\big)$ or $I_2-A\in M_2\big(P(R)\big)$, or $trA\in 1+P(R)$ and $tr^2A-4detA=u^2$ for a $u\in 1+P(R)$.}
\vspace{-.5mm}
\end{enumerate}\vspace{-.5mm} {\it Proof.}\ \ $(1)\Rightarrow (2)$ is clear by Theorem 5.1.

 $(2)\Rightarrow (1)$ If $trA\in 1+P(R)$ and $tr^2A-4detA=u^2$ for some $u\in 1+P(R)$,
 then $u\in U(R)$ and the equation $x^2-trA\cdot x+detA=0$ has a root
 $\frac{1}{2}(trA-u)$ in $P(R)$ and a root $\frac{1}{2}(trA+u)$ in
$1+P(R)$. Therefore we complete the proof by Theorem
4.4.\hfill$\Box$

\vskip4mm \no{\bf Example 5.3.}\ \ {\it Let $R$ be a commutative
local ring, and let $p\in P(R),q\in R$. Then $\left(
\begin{array}{cc}
p+1&p\\
q&p \end{array} \right)$ is strongly $P$-clean if and only if
$1+4pq=u^2$ for a $u\in 1+P(R)$.} \vskip2mm\hspace{-1.8em} {\it
Proof.}\ \ Set $A=\left(
\begin{array}{cc}
p+1&p\\
q&p \end{array} \right)$. Then $A, I_2-A\not\in
M_2\big(P(R)\big)$. As $tr^2A-4detA=1+4pq$, the result follows by
Theorem 5.1.\hfill$\Box$

\vskip4mm \no {\bf Theorem 5.4.}\ \ {\it Let $R$ be a commutative
local ring, and let $A\in M_2(R)$. Then $A$ is strongly $P$-clean
if and only if} \vspace{-.5mm}
\begin{enumerate}
\item [(1)]{\it $A\in M_2\big(P(R)\big)$, or} \vspace{-.5mm}
\item [(2)]{\it  $I_2-A\in M_2\big(P(R)\big)$, or}
\item [(3)]{\it $A\in M_2(R)$ is strongly $\pi$-regular and $A$ is similar
to a matrix $\left(
\begin{array}{cc}
0&\lambda\\
1&\mu\end{array} \right)$, where $\lambda\in P(R),\mu\in 1+P(R)$.}
\end{enumerate}\vspace{-.5mm} {\it Proof.}\ \ Let $A\in M_2(R)$ be strongly $P$-regular.
Assume that $A, I_2-A\not\in M_2\big(P(R)\big)$. In view of Lemma
4.3, there exists a $P\in GL_2(R)$ such that $P^{-1}AP=\left(
\begin{array}{cc}
0&\lambda\\
1&\mu\end{array} \right)$, where $\lambda\in P(R),\mu\in 1+P(R)$.
According to Theorem 4.4, the equation $x^2-trA\cdot x+detA=0$ has
a root in $P(R)$ and a root in $1+P(R)$. As $trA=\mu$ and
$detA=-\lambda$, we see that $h(x)=x^2-\mu x-\lambda$ has two left
roots, one is in $U(R)$ and one which is nilpotent. In light of
[11, Lemma 20], we conclude that $P^{-1}AP$ is strongly
$\pi$-regular. Thus, we can find some $m\in {\Bbb N}$ such that
$\big(P^{-1}AP\big)^{m}=\big(P^{-1}AP\big)^{m+1}B$ and
$(P^{-1}AP)B=B(P^{-1}AP)$. It follows that $A^m=A^{m+1}(PBP^{-1})$
and $A(PBP^{-1})=(PBP^{-1})A$, and thus $A\in M_2(R)$ is strongly
$\pi$-regular.

Conversely, assume that $A\in M_2(R)$ is strongly $\pi$-regular and $A$ is similar
to a matrix $\left(
\begin{array}{cc}
0&\lambda\\
1&\mu\end{array} \right)$, where $\lambda\in P(R),\mu\in 1+P(R)$.
Then $\left(
\begin{array}{cc}
0&\lambda\\
1&\mu\end{array} \right)$ is strongly $\pi$-regular where
$\lambda\in R$ is nilpotent and $\mu\in U(R)$. In light of [11,
Lemma 20], $x^2-\mu x-\lambda$ has two roots, one $\alpha\in U(R)$
and one $\beta\in R$ which is nilpotent. Obviously, $\alpha^2-\mu
\alpha-\lambda=0$ and $\beta^2-\mu \beta-\lambda=0$; hence,
$\alpha+\beta=\mu$. As $R$ is commutative, we see that $\beta\in
P(R)$, and then $\alpha=\mu-\beta\in 1+P(R)$. Obviously, $trA=\mu$
and $detA=-\lambda$. Therefore the equation $x^2-trA\cdot
x+detA=0$ has two roots, one in $1+P(R)$ and one which is in
$P(R)$. According to Theorem 4.4, $A$ is strongly
$P$-clean.\hfill$\Box$

\vskip4mm \no {\bf Proposition 5.6.}\ \ {\it Let $R$ be a
commutative ring, and let $A\in M_2(R)$. If $R/J(R)\cong {\Bbb
Z}_2$ and $J(R)$ is nilpotent, then $A$ is strongly $\pi$-regular
if and only if $A\in GL_2(R)$ or $A$ is nilpotent, or $A$ is
strongly $P$-clean.} \vskip2mm \no {\it Proof.}\ \ If $A\in
GL_2(R)$  or $A$ is nilpotent, then $A$ is strongly $\pi$-regular.
If $A$ is strongly $P$-clean, it follows from Corollary 5.5 that
$A$ is strongly $\pi$-regular. Conversely, assume that $A$ is
strongly $\pi$-regular, $A\not\in GL_2(R)$ and $A\in M_2(R)$ is
not nilpotent. As $J(R)=P(R)$, we see that $A\not\in
M_2\big(J(R)\big)$. By virtue of [11, Lemma 19], $A$ is isomorphic
to $\left(
\begin{array}{cc}
0&\lambda\\
1&\mu\end{array} \right)$, where $\lambda\in P(R),\mu\in R$. If
$\mu\in 1+P(R)$, it follows from Theorem 5.4 that $A\in M_2(R)$ is
strongly $P$-clean. If $\mu\in P(R)$, then $A^2$ is isomorphic to
$\left(
\begin{array}{cc}
\lambda&\lambda\mu\\
\mu&\mu+\mu^2\end{array} \right)$. This implies that $A^2\in
M_2\big(P(R)\big)$. Hence, $A\in M_2(R)$ is nilpotent, a
contradiction. Therefore the result follows.\hfill$\Box$

\vskip4mm \no {\bf Example 5.7.}\ \ {\it Let $A\in M_2\big({\Bbb
Z}_{2^n}[i]\big) (n\geq 1)$. Then $A$ is strongly $\pi$-regular if
and only if $A\in GL_2\big({\Bbb Z}_{2^n}[i]\big)$ or $A$ is
nilpotent, or $A$ is strongly $P$-clean.} \vskip2mm \no {\it
Proof.}\ \ Clearly, $J\big({\Bbb Z}_{2^n}[i]\big)=\big(1+i\big)$,
and that ${\Bbb Z}_{2^n}[i]/J\big({\Bbb Z}_{2^n}[i]\big)\cong
{\Bbb Z}_2$. Thus, ${\Bbb Z}_{2^n}[i]$ is a commutative local ring
with the nilpotent Jacobson radical. Therefore we complete the
proof by Proposition 5.6.\hfill$\Box$

\vskip15mm \bc {\bf Acknowledgements}\ec \vskip4mm \no This
research was supported by the Nature Science Foundation of
Zhejiang Province (LY13A010019) and the Scientific and
Technological Research Council of Turkey (2221 Visiting Scientists
Fellowship Programme).

\vskip15mm \bc{\bf REFERENCES}\ec \vskip 4mm \small{\re{1} G.
Borooah; A.J. Diesl and T.J. Dorsey, Strongly clean triangular
matrix rings over local rings, {\it J. Algebra}, {\bf 312}(2007),
773--797.

\re{2} G. Borooah; A.J. Diesl and T.J. Dorsey, Strongly clean matrix
rings over commutative local rings, {\it J. Pure Appl. Algebra},
{\bf 212}(2008), 281--296.

\re{3} J. Chen; X. Yang and Y. Zhou, When is the $2\times 2$ matrix
ring over a commutative local ring strongly clean? {\it J. Algebra},
{\bf 301}(2006), 280--293.

\re{4} J. Chen; X. Yang and Y. Zhou, On strongly clean matrix and
triangular matrix rings, {\it Comm. Algebra}, {\bf 34}(2006),
3659--3674.

\re{5} H. Chen, {\it Rings Related Stable Range Conditions},
Series in Algebra 11, World Scientific, Hackensack, NJ, 2011.

\re{6} H. Chen, On strongly nil clean matrices, {\it Comm.
Algebra}, 41:3(2013), 1074-1086.

\re{7} H. Chen; A. Harmanci and A.C. Ozcan, Strongly nil-*-clean
rings, arXiv:1211.5286v1, 2012.

\re{8} A.J. Diesl, Nil clean rings,
 {\it J. Algebra}, {\bf 383}(2013), 197--211.

\re{9} L. Fan and X. Yang, A note on strongly clean matrix rings,
{\it Comm. Algebra}, {\bf 38}(2010), 799--806.

\re{10} J.E. Humphreys, {\it Introduction to Lie Algebra and
Representation Theory}, Springer-Verlag, Beijing, 2006.

\re{11} Y. Li, Strongly clean matrix rings over local rings, {\it
J. Algebra}, {\bf 312}(2007), 397--404.

\re{12} W.K. Nicholson and Y. Zhou, Rings in which elements are
uniquely the sum of an idempotent and a unit, {\it Glasgow Math.
J.}, {\bf 46}(2004), 227--236.

\re{13} X. Yang and Y. Zhou, Strongly cleanness of the $2\times 2$
matrix ring over a general local ring, {\it J. Algebra}, {\bf
320}(2008), 2280-2290.

\end{document}